\documentclass[12pt]{amsart}

\usepackage{latexsym}
\usepackage{amssymb}
\usepackage{amsmath}

\newtheorem{theorem}{Theorem}[section]
\newtheorem{lemma}[theorem]{Lemma}
\newtheorem{corollary}[theorem]{Corollary}
\newtheorem{proposition}[theorem]{Proposition}
\newtheorem{remark}[theorem]{Remark}
\newtheorem{remarks}[theorem]{Remarks}
\newtheorem{definition}[theorem]{Definition}

\newcommand{\proofcl}{\noindent {\it {Proof of the Claim. }}}
\newcommand{\prend}{ $\Box $\hfill \bigskip}
\newcommand\spec{\mathop{\rm sp}}
\newcommand{\cl}[1]{\mathcal{#1}}
\newcommand{\bb}[1]{\mathbb{#1}}

\def\gs{\sigma}
\def\wee{\wedge}
\def\gl{\lambda}

\def\bh{\mathcal{B}(H)}

\newcommand{\nor}[1]{\left\Vert #1\right\Vert}

\title[Angles in C*-algebras]%
 {Angles in C*-algebras}

\author{M. Anoussis, A. Katavolos \and I. G. Todorov}

\address{M. Anoussis\\ Department of Mathematics,\\
University of the Aegean,\\
GR-83 200, Karlovassi,\\ Samos, GREECE} \email{mano@aegean.gr}
\address{A. Katavolos\\ Department of Mathematics,\\
University of Athens, \\
Panepistimioupolis, \\GR-157 84, Athens, GREECE}
   \email{akatavol@math.uoa.gr}
\address{I.G. Todorov\\ Department of Pure Mathematics,\\
Queen's University Belfast,\\
Belfast BT7 1NN,\\ UNITED KINGDOM}
   \email{i.todorov@qub.ac.uk}

\thanks{This project is cofunded by European Social Fund and National
Resources - (EPEAEK II) ``Pyhtagoras II'' grant No. 70/3/7997}


\begin{document}
\maketitle

\begin{abstract}
In this work we characterise the C*-algebras $\cl A$ generated by projections
with the property that every pair of projections in $\cl A$ has positive angle, 
as certain extensions of abelian algebras by algebras of compact operators.
We show that this property is equivalent to a lattice theoretic property of projections
and  also to the property that the set of finite-dimensional *-subalgebras of $\cl A$
is directed. 
\end{abstract}

\section{Introduction}

The structure of the set
of projections of an operator algebra has
been one of the main objects of investigation from the beginnings
of the theory.  In many cases, knowledge of the projections and
their equivalence classes provides useful information about the
structure of the algebra. The structure of projections in a C*-algebra is being
studied extensively. One can find a survey of results in
\cite{bla}. 

The C*-algebra generated by two projections
is analysed in Pedersen~\cite{ped} and in Raeburn and Sinclair~\cite{sin}.
The relative position of two  Hilbert space projections is studied 
in Dixmier~\cite{di}, Davis~\cite{davis} and Halmos~\cite{halm}. Some
lattice-theoretic implications of positive angles appear in
Mackey~\cite{mac}. Also,  consequences of this geometric
property for the structure of $\cl A$ may be found in Akemann~\cite{ak}. 
Actually, this paper motivated some of our results in Section 3.

D. Topping \cite{to} has characterised the von Neumann algebras
$\cl M$ with the property that every pair of projections in $\cl M$ has positive angle.
In this work we characterise the C*-algebras $\cl A$ 
generated by  projections and   having this property, 
which we call {\it the positive angle property}.

In Section 2 we define the notion of angle between two projections in
a C*-algebra $\cl A$ and show that it does not depend on the (faithful)
representation of $\cl A$.

We say that a C*-algebra $\cl A$ has the {\it sublattice property} if
the set $\cl{P(A)}$ of all projections in $\cl A$ is a sublattice
of the projection lattice of its enveloping von Neumann algebra
$\cl{A}^{**}$.  This property is studied in Section 3. The main result
is Theorem \ref{th_char2}
which asserts that the sublattice property is equivalent to the positive angle property.

In Section 4 we analyse the structure of a C*-algebra which satisfies
the positive angle property.  Our analysis has two main consequences:

We show that the positive angle property is equivalent
to the {\it directed set property}:
A C*-algebra $\cl A$ is said to have  this
property \cite{laztr} if the set of finite dimensional
C*-subalgebras of $\cl A$ is directed by inclusion.

Thus  a geometric,  a lattice-theoretic
and an algebraic property are shown to be equivalent.

The main result of our analysis 
concerns the structure of a C*-algebra $\cl A$ 
satisfying the positive angle property and generated by projections.
Such an algebra is characterized (Theorem \ref{thm}) as an extension \\
$0\to \cl K \to \cl A \to \cl C \to 0$
where $\cl K$ is a C*-algebra of compact operators, 
$\cl C$ is an abelian C*-algebra generated by projections, 
and the Busby invariant of the extension takes values in the 
centre of the corona of $\cl K$. Alternatively, $\cl A$ can be written
as a sum $\cl K+\cl Z$, where $\cl Z$ is its centre.   
 
In particular, we characterise the AF algebras having the positive angle
property, or equivalently the directed set property.
Lazar in \cite{laztr} has given a different characterisation  
of these algebras in spectral terms.

As a corollary, we show that Topping's characterisation of von Neumann algebras
with the positive angle property extends to AW* algebras.

\vspace{2ex}

\noindent\textbf{Notation} In general we follow the notation of
\cite{ped2}.  By $\bh$ we denote the algebra of all bounded linear
operators on a Hilbert space $H$ and by $\cl{K}(H)$ the
subalgebra of all compact operators on $H$.

If $\cl A$ is a C*-algebra we denote by $\cl{P(A)}$ the partially
ordered set  of all projections in $\cl A$. When $\cl A$ is a von
Neumann algebra (or more generally, an AW*-algebra), then
$\cl{P(A)}$ is a (complete) lattice. If  $p,q$ are projections, we
denote by $p\vee q$ their supremum and by $p\wee q$ their infimum;
for a family $\cl E$ we use the symbols $\vee \cl E$ and $\wee \cl
E$.  The partially ordered set where these are calculated will be
clear form the context.

The universal representation of $\cl A$ is denoted by $(\pi_u, H_u)$.
We often identify $\cl A$ with $\pi_u(\cl A)$, and the
second dual $\cl A^{**}$ with the bicommutant  $\pi_u(\cl A)''$.
If $X$ is a subset of a vector space we denote by $[X]$ the
linear span of $X$. For a family $\{X_j\}$ of Banach spaces, we use the symbol
$\sum_{j\in \cl J} \oplus_{\ell^\infty}X_j$ for the space of all bounded families 
$(x_j)$ with $x_j\in X_j$ and the symbol $\sum_{j\in \cl J} \oplus_{c_o}X_j$ 
for the subspace of families $(x_j)$ with $(\nor{x_j})\in c_o(\cl J)$.   

\section{Angles}

In this section we show that the notion of angle between
two projections in a C*-algebra can be defined intrinsically;
in particular, it is the same in any faithful representation
of the algebra (see Remark \ref{ivan}).

If $p,q$ are in a C*-algebra $\cl A$,  we
denote by $ C^*(p,q)$ the C*-subalgebra of $\mathcal{A}$ generated by $p$ and
$q$. We will use the following description of the C*-algebra
generated by two projections due to Pedersen \cite{ped}.

\begin{proposition}\label{peder}
Let $p,q$ be projections on a Hilbert space $H$ and let $\gs
=\spec (pqp), \, \gs_0=\gs\setminus \{0\}, \, \gs_1=\gs\setminus
\{1\}, \, \gs_{01}=\gs\setminus \{0,1\}$. Then the C*-algebra $\cl
C =C^*(p,q)$ is isomorphic to a subalgebra of $M_2(C(\gs))$. More
precisely, if $\cl{C}_o$ is the C*-subalgebra of $\cl C$ generated
by  $pq$ then
\begin{eqnarray*}
\cl C_o & \simeq & \left[
\begin{array}
[c]{cc}%
C_{o}(\gs_0) & C_{o}(\gs_{01})\\
C_{o}(\gs_{01}) & C_{o}(\gs_{01})
\end{array}
\right]  \\ &
 = & \left\{  \left[
\begin{array}
[c]{cc}%
f_{11} & f_{12}\\
f_{21} & f_{22}%
\end{array}
\right]  :f_{11}(0)=0, f_{ij}(0)=0=f_{ij}(1) \mbox{ for
}ij>1\right\} .
\end{eqnarray*}
(a) If $\ 0\in\overline{\gs_0}$, then
\begin{eqnarray*}
\mathcal{C}  & \simeq & \left[
\begin{array}
[c]{cc}%
C(\gs) & C_{o}(\gs_{01})\\
C_{o}(\gs_{01}) & C_{o}(\gs_{1})
\end{array}
\right]  \\
& = & \left\{  \left[
\begin{array}
[c]{cc}%
f_{11} & f_{12}\\
f_{21} & f_{22}%
\end{array}
\right]  :f_{ij}(0)=0=f_{ij}(1)\mbox{ for }i\neq j,
f_{22}(1)=0\right\}  .
\end{eqnarray*}
(b) If $0\notin\overline{\gs_0}$, then
$$\cl{C}=\cl{C}_{o}\oplus\mathbb{C} (p\wedge
q^{\bot})\oplus\mathbb{C}(p^\bot \wedge q).$$
\end{proposition}

The next proposition will be used
repeatedly throughout the paper.

\begin{proposition}\label{isod}
Let $p,q\in\bh$ be projections, $e=p\wedge q$ and  $f=p\vee q$. 
The following are equivalent:

(a) The sequence $\{(pqp)^{n}\}_n$ is norm-convergent (to $e$);

(b) The point $1$ is not an accumulation point of $\spec (pqp)$;

(c) $\left\|  (p-e)(q-e)\right\|  <1$;

(d) $f\in C^*(p,q)$;

(e) $C^*(p,q)$ is unital.
\end{proposition}

\proof
(a)$\Rightarrow$(b) It is well known that $p\wedge q$ is
the strong limit of the sequence $\{(pqp)^{n}\}_n$ (in
$\bh$). Thus if $(pqp)^{n}$ is norm-convergent, it
converges to $e$.

Now $pqp$ is a positive contraction,
hence $C^*(pqp)\simeq C_{o}(\spec(pqp)\setminus\{0\})$.

If $\phi\in C_{o}(\spec (pqp)\setminus\{0\})$ denotes the function
$\phi(t)=t$, then $\phi^{n}$ decreases pointwise to $\chi_{\{1\}}$
(the characteristic function of ${\{1\}}$). Thus if $(pqp)^{n}$ converges in norm,  
then $\phi^n$ converges uniformly and so
$\chi_{\{1\}}$ must be continuous on $\spec (pqp)$; hence $1$
cannot be an accumulation point of $\spec (pqp)$.

\medskip

(b)$\Rightarrow$(c) If $1$ is not an accumulation point of
$\mathop{\rm sp} (pqp)$, then, in the notation of the previous paragraph, $\nor{(1-\chi_{\{1\}})\phi}_\infty
<1$; 
thus $\|(1-e)pqp\|<1$ and so
\begin{align*}
\|(q-e)(p-e)\|^2 & = \|(p-e)(q-e)(p-e)\| =\|pqp-e\| \\
 & = \|pqp(1-e)\|<1.
\end{align*}
(c) $\Rightarrow$ (a) If $\|(q-e)(p-e)\|<1$ then, as just observed,
$(p-e)(q-e)(p-e)$ is a strict contraction and so
$((p-e)(q-e)(p-e))^{n}$ is
norm-convergent to $0$; but it is easily verified that 
$((p-e)(q-e)(p-e))^n=(pqp)^n-e.$

\medskip

(b)$\Rightarrow$(e) If $1$ is not a limit point of $\spec(pqp)$, then   the
characteristic function $\chi$ of $\sigma_1$ is in $C_0(\sigma_1)$. 
Thus if $0$ is  a limit point of $\sigma(pqp)$, then from part (a) of \ref{peder},
the element  $$u=\left(\smallmatrix 1 & 0\\
0 & \chi\endsmallmatrix\right)\, ,$$
is in the algebra and  clearly acts as a unit.
 If $0$ is not a limit point of $\sigma(pqp)$, then from part (b) of \ref{peder},
the element  $$u=\left(\smallmatrix 1 & 0\\
0 & \psi\endsmallmatrix\right)\oplus p\wedge q^\bot\oplus p^\bot\wedge q\, ,$$  
where $\psi$ is the characteristic function of $\sigma_{0,1}$ is in the algebra 
and acts as a unit for the algebra.

\medskip

(e)$\Rightarrow$(b) We use Proposition \ref{peder}:
Let $u=[u_{ij}]\in M_2(C(\gs))$ denote the unit of
$\mathcal{C}\simeq C^{*}(p,q)$. If $g=[g_{ij}]\in \cl C$  is
chosen so that $g_{12}=0$, then the relation $ug=g$ yields
$u_{22}g_{22}=g_{22}$. Since $g_{22}\in C_o(\gs_1)$ is arbitrary,
it follows that $u_{22}(t)=1$ for all $t\in \gs \setminus \{1\}$.
Since $u_{22}\in C_o(\gs_1)$, the point $1$ cannot be an accumulation point of
$\gs=\spec (pqp)$.

\medskip

(e)$\Leftrightarrow$(d)  If $f\in C^*(p,q)$, then  since $fp=pf=p$
and $fq=qf=q$ it is clear that $f$ is  the unit of $C^*(p,q)$.

Conversely, let $u$ be the unit of $C^*(p,q)$. Then $u\ge p$ and $u\ge q$,
so $u\ge f$ (in $\bh$).
Now $u$  can be approximated in norm by polynomials $\phi_n$ in $p$ and
$q$. But since $fp=pf=p$ and $fq=qf=q$ we have $\phi_nf=f\phi_n = \phi_n$
for each $n$, so in
the limit $fu=uf=u$. Thus $f\ge u$ and so $f=u$.

The proof is complete. \prend

\begin{remark}
The conditions of the proposition are also equivalent to

(i) $\left\|  (f-p)(f-q)\right\|  <1$ \\
and to

(ii) The sequence $\{(p+q)^{1/n}\}_n$ is norm-convergent (to $f$).\\
Since we will not use this fact, we omit the proof.
\end{remark}

If $\pi$ is a representation of a c*-algebra $\cl A$, we denote by
$\tilde{\pi}$ its unique extension  to a normal (i.e.
w*-continuous) representation of $\cl A^{**}$.

\begin{proposition}\label{pr_ang} Let $p,q$ be projections in a C*-algebra $\cl{A}$ and
$e=p\wedge q$, calculated in $\cl A^{**}$.
If  $\pi$ is a faithful representation of $\cl A$ on a Hilbert
space $H$  then

(i)  $\tilde{\pi}(e) \;\mbox{is the infimum of $\pi(p)$ and
$\pi(q)$ calculated in }\; \bh $;

(ii)
$\Vert(\pi(p)-\tilde{\pi}(e))(\pi(q)-\tilde{\pi}(e))\Vert=\|(p-e)(q-e)\|
$.

\end{proposition}
\begin{proof}
Let $e_1= \pi(p)\wedge\pi(q)$. 
Since $e_1 = \lim_n (\pi(p)\pi(q)\pi(p))^n$ and $e=\lim_n (pqp)^n$ strongly,
(i) follows from the normality of $\tilde\pi$. 

Hence
$$
\Vert(\pi(p)-e_1)(\pi(q)-e_1)\Vert=\|\tilde{\pi}((p-e)(q-e))\| \le
\|(p-e)(q-e)\|\le 1. $$

If $\Vert(\pi(p)-e_1)(\pi(q)-e_1)\Vert=1$ then we have equality
throughout.

If on the other hand $\Vert(\pi(p)-e_1)(\pi(q)-e_1)\Vert<1$ then, by
Proposition \ref{isod},  $e_1=\tilde{\pi}(e)$ lies in the C*-algebra
$\pi(\cl A)$. Thus there exists $r\in \cl A$ with $\tilde{\pi} (e)
= \pi (r)$ and so
$$\|\pi(r) -  (\pi(p)\pi(q)\pi(p))^n\| \to 0 .$$
It follows that $\|r - (pqp)^n\| \to 0.$ Since $e$ is the
strong operator limit of $(pqp)^n$, it must equal $r$; thus $e\in
\cl A$. But $\pi$ is isometric, so $$
\Vert(\pi(p)-\tilde{\pi}(e))(\pi(q)-\tilde{\pi}(e))\Vert=
\Vert(\pi(p)-\pi(e))(\pi(q)-\pi(e))\Vert=\|(p-e)(q-e)\| $$ as
required. \end{proof}

\begin{definition}
If $\cl A$ is a C* algebra and $p,q \in \cl {P(A)}$, we define
 $$c(p,q) = \|(p-p\wedge q)(q-p\wedge q)\|$$
where $p\wee q$ is calculated in $\cl A^{**}$.
\end{definition}

\begin{remarks}\label{ivan}
\emph{ \textbf{(i)} Let $p,q\in\bh$ be two  projections with
$pH \cap qH=0$. Recall that the angle $\theta(p,q)$ between the
subspaces $pH$ and $qH$ is defined to be the arc cosine of
the quantity $$ \sup\{|(\xi, \eta)| :  \xi = p\xi,  \eta = q\eta,
\|\xi\|,\|\eta\| \leq 1\}.$$ Note that $$\cos\theta(p,q) = \|pq\|.$$ It
is well known that $\theta(p,q) > 0$ if and only if $pH  +
qH $ is closed.}

\smallskip

 \noindent  \emph{\textbf{(ii)} It follows from Proposition \ref{pr_ang}
that if $p$ and $q$ are projections in a C*-algebra $\cl A$,
$$c(p,q) =\|(\pi(p)-\pi(p)\wedge \pi(q))(\pi(q)-\pi(p)\wedge
\pi(q))\|$$ for any faithful representation $\pi$ of $\cl A$. In
particular, if $\cl{A}$ is contained in a C*-algebra $\cl B$, then
$c(p,q)$ is independent of whether $p\wedge q$ is calculated in
$\cl A^{**}$ or $\cl B^{**}$.}

\smallskip

\noindent \emph{\textbf{(iii)} In Proposition \ref{isod} it is
shown that  $c(p,q)=1$ if and only if $1$ is an accumulation point
of $\spec (pqp)$. Applying the Proposition to the pair
$(p,q^{\perp})$ we see that $c(p,q^\bot)=1$ if and only if $0$ is
an accumulation point of $\spec(pqp)$.}
\end{remarks}

\section{Lattices of projections}

A C*-algebra $\cl A$ is said to have  the {\it directed set
property} \cite{laztr} if the set of finite dimensional
C*-subalgebras of $\cl A$ is directed by inclusion. It is said to
have the {\it lattice property} \cite{lazsc} if the partially
ordered set $\cl{P(A})$ of projections in $\cl A$ is a lattice in
its own order.  A. Lazar
\cite{lazsc} has shown that these properties are equivalent for AF
C*-algebras. In general, the lattice property does not imply the
directed set property; consider for example the C*-algebra $\bh$.

We introduce a lattice-theoretic property which will prove to be equivalent
to the directed set property:  We say that
$\cl A$ has the {\it sublattice property} when $\cl{P(A})$ is a
sublattice of $\cl{P(A}^{**})$ (of course $\cl{P(A}^{**})$ is  always a
lattice, since $\cl A^{**}$ is a von Neumann algebra).

The main result of this section is Theorem \ref{th_char2} which
shows that a C*-algebra $\cl A$ satisfies the sublattice property if and only if
$c(p,q) <1$ for each pair of projections $p,q\in \cl A$; this
gives a geometric characterization of the algebras with the
sublattice property. In Section 5 we will prove that this latter property is
equivalent to the directed set property (Theorem \ref{dsp}).

If $\pi$ is a faithful representation of a C*-algebra $\cl A$ and $p,q\in\cl A$
are projections such that $\pi(p)$ and $\pi(q)$ satisfy the equivalent
conditions of Proposition \ref{isod}, then $\pi(p)\wee\pi(q)\in C^*(\pi(p),\pi(q))$.
The converse does not hold in general; however, it is true if $\pi$ is the universal representation:

\begin{proposition}\label{p_spanclosed}
Let $\cl A$ be a C*-algebra,  $p,q\in \cl A$ be projections and
$f=p\vee q, \, e=p\wee q$ calculated in $\cl A^{**}$.
 The following are equivalent:\label{rem33} \\
 (i) $e\in \cl A$ \\
(ii) $c(p,q) < 1$\\
(iii) $f\in\cl A$.
\end{proposition}
\begin{proof}  
Proposition \ref{isod} shows that $(iii)$ is equivalent to $(ii)$ 
and $(ii)$ implies $(i)$.  

$(i) \Rightarrow (ii)$ Assume that $c(p,q)=1$ and that $e\in\cl
A$. Then the projections $p_o=p-e$ and $q_o=q-e$ are in $\cl A$.
Now $c(p_o,q_o)=1$ and hence 
$1$ is an accumulation point of
$\spec (p_oq_op_o)$, by Proposition \ref{isod}. The C*-algebra generated by
$\{p_oq_op_o,p_o\}$ consists of continuous functions on $\spec
(p_oq_op_o)$; let $\varphi$ be the state of this algebra
corresponding to evaluation at $1$. Now define a state $\psi$ of the
C*-algebra generated by $\{p_o,q_o\}$ by
$\psi(x)=\varphi(p_oxp_o)$, and extend it to a state 
$\tilde{\psi}$ of $\cl A$. By hypothesis, $\cl A$ is sitting in
its universal representation, hence $\tilde{\psi}$ is a vector
state. Since $\tilde{\psi}(p_o)=\tilde{\psi}(q_o)=1$, it follows
that $\tilde{\psi}(p_0\wee q_o)=1$; this is impossible, since
$p_o\wee q_o=0.$ 
 \end{proof}

\bigskip

Let $\cl A$ be a C*-algebra. Recall \cite{ak} that a projection
$p\in \cl A^{**}$ is said to be \textit{open} if there exists an
increasing net ${a_i}\in \cl A^+$ such that $a_i\nearrow p$. A
projection $p\in \cl A^{**}$ is said to be \textit{closed} if
$1-p$ is open. In case $\cl A$ is unital, this means that $p$ is
the infimum of a decreasing net of contractions from $\cl A$.

Akemann \cite[Theorem II.7]{ak} proves that if $p,q$ are closed
projections in $\cl A^{**}$ and $c(p,q)<1$, then $p\vee q$ is
closed. The converse is not true in general, as we show below.
However,  the following is immediate from  Proposition
\ref{p_spanclosed}.

\begin{corollary}
Let $\cl A$ be a unital C*-algebra,  $p,q\in \cl A$ be projections
and $f=p\vee q$ calculated in $\cl A^{**}$. The projection $f$ is
closed if and only if $c(p,q)<1.$
\end{corollary}
\begin{proof}
The supremum of a set of projections in $\cl A$ is always
open \cite[Proposition II.5]{ak}. So, if $f$ is closed then
it must belong to $\cl A$  \cite[Proposition II.18]{ak}.
By Proposition \ref{p_spanclosed}, this implies $c(p,q)<1$.
\end{proof}

The following proposition shows that this Corollary does not hold
if we do not assume that $p, q \in \cl A$.
\begin{proposition}
There exists a unital C*-algebra $\cl A$ and open projections $p$
and $q$ in $\cl A^{**}$
 such that $c(p,q)=1$ and
$p\wedge q$ is open.
\end{proposition}
\begin{proof} Let $\cl A=\mathbb{C}I+\cl{K}(H)$. If
$\{\xi_n\}_{n\in\bb N}$ is orthonormal in $H$, consider the
following compact operators
$$a=\sum_{n=1}^{\infty}\frac{1}{n}e_{n},\quad\qquad b=\sum_{n=1}^{\infty}\frac{1}{n}f_{n},$$
where $e_{n}, \ f_{n}$ are the rank one projections with ranges
spanned by $\xi_{2n}$ and $\eta_{n}=\cos (\frac 1n) \xi_{2n}+\sin
(\frac 1n)\xi_{2n+1}$, respectively.
Let $p$ be  the strong limit in $\cl A^{**}$ of the sequence
$a^{1/k}$ and $q$ the strong limit  of the sequence $b^{1/k}$.
Then  $p$ and $q$ are open projections by construction and
it is clear that $c(p,q)=1$. We claim that $p\wedge q=0$ and hence
$p\wedge q$ is open. Indeed $p$ and $q$ are the range projections
of $\pi_u(a)$ and $\pi_u(b)$ respectively. These operators are
ampliations of $a$ and $b$ \cite[Corollary 1 of Theorem
1.4.4]{arv} hence the closures of their ranges do not intersect. \end{proof}

The following lemma will be needed in the proof of Theorem
\ref{th_char2}.

\begin{lemma}\label{l_sp}
Let $p$ and $q$ be projections on a Hilbert space $H$ and $\cl
A$ be the C*-algebra generated by $p$ and $q$. If $\spec(pqp)$ is
infinite then $\cl A$ contains two distinct 
projections $r$ and $s$ with $c(s,r)=1$.
\end{lemma}
\begin{proof} If $c(p,q)=1$  we are done. Assume that $c(p,q) <1$.
This implies that  $\cl A$ is unital and that $1$ is not an
accumulation point of $\spec (pqp)$ (Proposition \ref{isod}). If
$0$ is an accumulation point of $\spec (pqp)$, then $c(p^\bot ,
q)=1$ and again we are done. Thus we may assume that $K=\spec
(pqp)\setminus \{0,1\}$ is a compact subset of $(0,1)$.

Since $\spec (pqp)$ is infinite, it has an accumulation point
$\delta\in K$. Using the description of the algebra $\cl A$ given
in Pedersen's theorem (Proposition \ref{peder}), we construct
projections $r,s\in \cl A$ with $c(s,r)=1$.
Since $0$ is not an accumulation point of $\gs=\spec(pqp)$,
the characteristic function $\chi$ of $\gs_{0}$ is in $C_0(\gs_0)$; thus
by part (b) of Pedersen's theorem, the projection
$s=\left(\smallmatrix \chi & 0\\
0 & 0\endsmallmatrix\right)$ is in $\cl A$.

Choose a continuous function $f:\spec (pqp)\to [0,1]$ supported in
$K$, with $f(\delta) = 1$ and such that $f(\delta)$ is an
accumulation point of $f(K)$. Define $g$ on $\spec (pqp)$ by
$g(t)=(1-f^2(t))^{\frac 12}$ for $t\in K$ and $g(t)=0$ for $t\in
\spec (pqp)\setminus K$. Now consider
$r = \left(\smallmatrix f^2 & fg\\
fg & g^2\endsmallmatrix\right)$. This is a projection in $\cl A$.
Note that $srs = \left(\smallmatrix f^2 & 0\\ 0 & 0 \endsmallmatrix\right)$.
It follows that 1 is an accumulation point of $\spec(srs)$ and so $c(s,r) = 1$.
\end{proof}

\begin{remark}\label{rfd}
Consider the C*-algebra $C^*(p,q)$ generated by two projections. 
If $c(r,s)<1$ for all projections $r,s\in C^*(p,q)$, then the algebra 
is finite dimensional. Indeed, $\spec (pqp)$ is
finite by Lemma \ref{l_sp}, and the claim follows from Proposition \ref{peder}. 

We will prove in Theorem \ref{52} that the same holds for the C*-algebra
generated by any finite number of projections.  
\end{remark}

\begin{theorem}\label{th_char2}
Let $\cl A$ be a C*-algebra. The following are equivalent:

(i) The partially ordered set  $\cl P(\cl A)$ is a sublattice of
the lattice $\cl P(\cl A^{**})$;

(ii) $c(p,q) <1$ for each pair $p,q\in \cl A$ of projections;

(iii) $\spec(pqp)$ is finite, for each pair $p,q\in \cl A$ of
projections.
\end{theorem}
\begin{proof} 
The implications $(ii) \Rightarrow (i)$ and $(iii) \Rightarrow (ii)$
follow from Proposition \ref{isod} and $(i) \Rightarrow (ii)$ 
follows from Proposition \ref{rem33}.

We prove $(ii) \Rightarrow (iii)$: If $\spec(pqp)$ is infinite 
for some projections $p,q\in\cl A$, then Lemma \ref{l_sp} shows
that there exist projections $p_1,q_1\in \cl A$ such that
$c(p_1,q_1) = 1$.
\end{proof}

\section{The Positive Angle Property}

We will say that a C*-algebra $\cl A$ has \emph{the positive angle property} when 
$c(p,q)<1$ for all $p,q\in \cl{P(A)}$. 
Two examples are immediate: abelian C*-algebras and the algebra of all compact operators,
or subalgebras of these.
Of course, projectionless C*-algebras have the property in a trivial way, so we will
be considering algebras  generated by their projections. 
We will show that the most general C*-algebra with these two properties 
can be constructed from an abelian C*-algebra and a C*-algebra of compacts
in the way described in Theorem \ref{thm}.

Thus in this section,   and up to Proposition \ref{gen}, 
\emph{we will let $\cl{A}$ 
be a C*-algebra with the positive angle property
which is generated, as a C*-algebra, by its projections.} 
 
Note that if $p,q$ are projections in $\cl A$ then 
$C^*(p,q)$ is finite-dimensional (Remark \ref{rfd}) 
and hence w*-closed in $\cl A^{**}$.

\medskip

It will be convenient to identify $\cl A$ with its reduced atomic representation acting
on a Hilbert space $H$. Thus if $\cl I=\widehat{\cl A}$ is the spectrum of $\cl A$,
and a representative $\pi_i$ is chosen from each unitary equivalence class $i\in\cl I$,
we have
\begin{align*}
\mathop{\rm id} = & \sum_{j\in \cl I} \oplus (\pi_j, H_j) \\
\text{and }\quad
\cl A'' = & \sum_{j\in\cl I}\oplus \cl B(H_j)
\end{align*}
(see for example \cite[10.3.10]{kr}.)

Denote by $z_j\in\cl B(H)$ the projection onto $H_j$. The $z_j$ 
are pairwise orthogonal central projections in
$\cl{A}''$. Also, $\pi_j(a)=az_j$ for each $a\in \cl A$ and $j\in\cl
I$.

Denote by $\cl{K(A)}$ the \emph{compact elements} of $\cl{A}$, that
is, the elements $a\in\cl A$ for which the operator $x\to axa$ is
compact on $\cl A$. It follows from \cite{yl} and \cite[Theorem
3.7]{erd} that $\cl K(\cl A) = \cl K(H)\cap\cl A$. We will prove in
Lemma \ref{l_twop} below that, if $\cl A$ is not abelian, the ideal $\cl{K(A)}$ is nonzero.
Let $\cl J = \{j\in\cl I:\pi_j(\cl{K(A)})\ne 0\}$. Note that
$j\in\cl J$ iff $\cl{K(A)}z_j\ne 0$.

From \cite[2.11.2]{di_cstar},  the spectrum of $\cl{K(A)}$ is in bijective
correspondence with $\cl J$.
It now follows by standard arguments (see eg. \cite[1.4.5]{arv}) 
that%
$$\cl K(\cl A) =\sum_{j\in \cl J} \oplus_{c_o}\cl{K(A)}z_j\oplus 0 =
\sum_{j\in \cl J} \oplus_{c_o}\cl K_j\oplus 0$$ where each $\cl
K_j=\cl{K(A)}z_j$ is an ideal of $\cl A$ and is equal to  $\cl
K(H_j)$.

\begin{remark}\label{compt}
 The centre $\cl{Z(A)}$ acts
`componentwise' on $\cl{K(A)}$: if $a\in \cl{Z(A)}$ and  $j\in \cl J$ there is a
scalar $\gl_j(a)$ such that $$az_j=\gl_j(a)z_j.$$
Also, the intersection $\cl{K(A)}\cap \cl{Z(A)}$ is isomorphic
to $c_o(\cl J_f)$ where
$\cl J_f=\{j\in \cl J : \dim H_j<\infty\}$.
\end{remark}
Indeed, since $\cl K_j\subseteq \cl A$ it follows that $a$ commutes with $\cl K_j$;
hence so does $az_j$. Since $\cl K_j=\cl K(H_j)$, $az_j$ is a scalar
multiple of the identity $z_j$ of $H_j$.

The second claim is obvious.
\begin{lemma}\label{l_twop}
If $f\in\cl A$ is a minimal projection then $f\cl A f = \bb{C}f$.
Thus, the minimal projections in $\cl A$ are precisely the rank
one projections of $\cl K_j$, $j\in\cl J$.
\end{lemma}
\proof  First observe that if $p,q\in \cl{P(A)}$,
then  $C^*(p,q)$  is finite dimensional, hence 
is linearly generated by its projections. 
Thus
$$pq=\sum_{i=1}^n\gl_ir_i, \;\; r_i\in \cl{P(A)}, \; \gl_i\in\bb C.$$
An obvious induction shows that any word $w=p_1p_2\dots p_n$ can be written 
\begin{align*}
w =  \sum_{i=1}^n\gl_ig_i, \;\;g_i\in \cl{P(A)}\quad 
\text{and so}\quad  fwf =  \sum_{i=1}^n\gl_ifg_if.
\end{align*}
Since each $C^*(f,g_i)$ is finite dimensional, by spectral theory $fg_if$ 
is a linear combination of projections in $\cl A$,
necessarily no larger than $f$. Since $f$ is minimal, each $fg_if$ is 
therefore a scalar multiple of $f$, hence so is $fwf$. 
Since $\cl A$ is generated by $\cl{P(A)}$, 
the first claim is proved.

Note that the minimal projections of $\cl K_j$ are clearly minimal
in $\cl A$. Conversely, by the previous paragraph, every minimal
projection $p$ of $\cl A$ is compact and hence belongs to some
$\cl K_j$. Clearly, $p$ must be of rank one. \prend

Recall that two projections $p,q\in\cl{P(A)}$ are said to be
(Murray - von Neumann) equivalent
if there exists a partial isometry $v\in\cl A$
such that $v^*v = p$ and $vv^* = q$; we write $p\sim q$.

\begin{proposition}\label{l_zeroa}
If $p\in \cl{P(A)}$, there does not exist an infinite orthogonal
family $\{e_n,
f_n\}_{n=1}^{\infty}\subseteq\cl A$ of nonzero 
projections in $\cl A$ such that $e_n\le p,\, f_n\le p^\perp$ and
$e_n\sim f_n$, for each $n\in\mathbb{N}$.
\end{proposition}
\proof Suppose such a family exists. We will reach a contradiction by
constructing a projection $q\in\cl A$ such that $c(p,q)=1$. 

The sum $\sum_kpe_k$ converges
strongly to a projection in $\cl A''$. Define $p_1=p-\sum_kpe_k$, so
that
$$p=\sum_kpe_k+p_1=\sum_ke_k+p_1$$ 
since $e_k\le p$. If $v_kv_k^*=e_k$
and $v_k^*v_k=f_k$, 
set $q_k=c_k^2e_k+c_ks_k(v_k+v_k^*)+s_k^2f_k$ where $c_k,s_k$ are real 
scalars with $c_k^2+s_k^2=1$ and $c_k\to 1$. The $q_k$ are
pairwise orthogonal projections and are all orthogonal to $p_1$
(because $p_1\le p$ while $f_k\le p^\perp$ so $p_1\bot f_k$, 
and $p_1\bot e_k$ by construction).
Define
$$q=\sum_kq_k+p_1.$$
We show that $q\in \cl A$. We have 
$$q=\sum_{k=1}^\infty q_k +p_1
= \left(\sum_{k=1}^\infty e_k +p_1\right)+ \sum_{k=1}^\infty (q_k-e_k)$$ and
the first term is in $\cl A$ (it is equal to $p$) while the second
converges in norm.

Indeed, it is in fact norm Cauchy, since
\begin{align*}
\nor{e_k-q_k} & = \nor{(1-c_k^2)e_k-c_ks_k(v_k+v_k^*)-s_k^2f_k} \\
& = s_k\nor{s_ke_k-c_k(v_k+v_k^*)-s_kf_k} \le 4s_k \to 0
\end{align*}
and
$$\nor{\sum_{k=m}^n (q_k-e_k)} = \max_{m\le k\le n}\nor{q_k-e_k}.$$

Setting $g_k=e_k+f_k$, it is easy to see that $g_kp = pg_k = e_k$
and $g_kq = qg_k = q_k$ for all $k$. Hence
$$(p\wedge q) g_k = g_k (p\wedge q) = (pg_k)\wedge (qg_k) =
e_k\wedge q_k = 0.$$ Thus, $p\wedge q = p_1$. It follows that
$c(e_k,q_k)=\|e_kq_k\|\leq \|pq-p_1\|=c(p,q)$ for each $k$ and so
$c(p,q) = 1$. \prend

\begin{lemma}\label{nc1}
If $p\in\mathcal{P(A)}$  is not central then there exist minimal
projections $e,f\in\cl{P(A)}$ such that $e\leq p$, $f\leq
p^{\perp}$ and $e\sim f$.
\end{lemma}
\proof First, we claim that there exist non-zero projections
$e_0,f_0\in\cl A$ such that $e_0\leq p$, $f_0\leq p^{\perp}$ and
$e_0\sim f_0$.

Indeed, since $\mathcal{P(A)}$ generates
$\mathcal{A}$, 
there exists $r\in\mathcal{P(A)}$ with $prp^{\bot}\neq0$. Let
$prp^{\bot}=v|prp^{\bot}|$ be the corresponding polar
decomposition. Then $v$ belongs to $C^*(p,r)$. Set
$e_0=vv^{\ast}\le p$ and $f_0=v^{\ast}v\le p^\perp$.

Now suppose that $e_0$ does not majorise a minimal projection.
Then there is an infinite strictly decreasing sequence $(q_n)$ of proper
nonzero subprojections of $e_0$; set $e_1=e_0-q_1$ and
$e_n=q_{n-1}-q_{n}$
for $n\ge 2$. If $f_n = v^*e_nv$, the family $\{e_n,f_n\}$ is orthogonal
and satisfies $e_n\le e_0, \, f_n\le e_0^\perp$ and $e_n\sim f_n$
for all $n$. By Proposition \ref{l_zeroa}, this is a contradiction.

Thus there exists a minimal subprojection $e$ of $e_0$; set $f
=v^*ev$. $\quad\Box$

\begin{proposition}\label{gen} Each   
projection in $\cl A$ is a linear combination of central
and finite rank projections. In fact 
$\cl{P(A)}\subseteq[\cl{P(K(A))}] + \cl{P(Z(A))}$.
\end{proposition}
\proof
Let $p\in\cl{P(A)}$. We claim that the set
$$\cl J_0 = \{j\in \cl J :  0 < pz_j < z_j\} 
$$ is finite.
Indeed, for each $j\in\cl J_0$, choose projections
 $e_j,f_j$ of rank one  with $e_j\le pz_{j}$ and
$f_j\le z_{j}- pz_{j}$. If $\cl J_0$ is infinite, we obtain an infinite family
$\{e_j,f_j:j\in\cl J_0\}\subseteq\cl A$ which violates 
 Proposition \ref{l_zeroa}.

Note that, for each $j$, the projections $pz_j$ and $p^\perp z_j$
cannot both have infinite rank. Indeed,  if they did, since
the C*-algebra $\cl A z_j$
contains $\cl K(H_j)$, 
there would be rank one projections $\{e_n,f_n, n\in\bb N\}$ in
$\cl A z_j$ with $e_n\le p$ and $f_n\le p^\perp$. As in Proposition \ref{l_zeroa},
this contradicts the positive  angle property.

Set $\cl J_1=\{j\in \cl J_0: \dim pH_j < \infty\}$ and
$\cl J_2= \cl J_0 \setminus \cl J_1$.
If $j \in \cl J_1$ (resp. $j\in\cl J_2$), then the
projection $pz_j$ (resp. $p^{\perp}z_j$) has finite rank hence is in $\cl A$.
Put $p_1=\sum_{j \in \cl J_1}pz_j$ and
$p_2=\sum_{j \in \cl J_2}p^{\perp}z_j$. These  are finite
orthogonal sums of finite rank projections in $\cl A$, hence belong to $\cl{K(A)}$.
Since $p-p_1$ and $p_2$ are  orthogonal projections, the sum
$q=p-p_1+p_2$ is in $\cl{P(A)}$.

We claim that
$q\in\cl{Z(A)}$. If not then Lemma \ref{nc1} yields the existence
of minimal equivalent projections $e,f\in\cl A$ such that $e\leq
q$ and $f\leq q^{\perp}$. By Lemma \ref{l_twop}, there exists $j$
such that $e\leq z_j$ and $f\leq z_j$. But by construction $qz_j$
is equal to either zero or $z_j$, a contradiction. $\qquad\Box$

\medskip
In what follows we will need the notion of the \textit{Busby invariant}, 
introduced in \cite{bus}.
Recall (see for instance \cite[Section 3.2]{wo})
that if $\cl K$ is an ideal of a C* algebra $\cl A$ and
$\cl C$ is the quotient $\cl A/\cl K$,
the \textit{Busby invariant} of the extension $0\to\cl K\to\cl A\to\cl C\to 0$
is the morphism $\tau: \cl C\to M(\cl K)/\cl K$ defined as follows: Let $\gs:\cl A\to M(\cl K)$
be the unique extension of the inclusion map $\cl K\to M(\cl K)$. Given $c\in\cl C$
let $a\in\cl A$ be any lift of $c$; then one sets $\tau(c)=\pi(\gs(a))$, where
$\pi:M(\cl K)\to M(\cl K)/\cl K$ is the quotient map. 
We will use this notation below.

\begin{theorem}\label{thm}
Let $\cl A$ be a C*-algebra. The following are equivalent:

(i) \ $\cl A$ has the positive angle property and is  generated by its projections;

(ii) \ $\cl A=\cl{K(A)}+ \cl{Z(A)}$ and $\cl{Z(A)}$
is generated by its projections;

(iii) 
$\cl A$ is an extension of $\cl C$ by $\cl K$, 
where $\cl C$ is an abelian C*-algebra  
generated by its projections,  $\cl K$ is a C*-algebra of compact
operators, and the Busby invariant 
$\tau : \cl C\rightarrow M(\cl K)/\cl K$ takes values 
in the centre of $M(\cl K)/\cl K$.
\end{theorem}

\proof 
$(i)\Rightarrow (ii)$ If $p\in\cl{P(A)}$ then $p\in [\cl P(\cl K(\cl A))] + \cl{P(Z(A))}$
(Proposition \ref{gen}). Since $\cl{K(A)} + C^*(\cl{P(Z(A))})$ is
a C*-subalgebra of $\cl A$ and the latter 
is generated  by its projections, we have
$$\cl A \subseteq \cl K(\cl A) + C^*(\cl{P(Z(A))}) \subseteq
\cl K(\cl A) + \cl Z(\cl A) \subseteq\cl A $$
and hence equality holds.
Since
$\cl{A}=\cl{K(A)}+ C^*(\cl{P(Z(A))})$, it follows that
$$\cl{Z(A)}=(\cl{Z(A)}\cap\cl{K(A)})+ C^*(\cl{P(Z(A))}).$$
But
$\cl{Z(A)}\cap\cl{K(A)}\simeq c_o(\cl J_f)$ (Remark \ref{compt})
and so $\cl{Z(A)}$ is generated by its projections.

\medskip

$(ii) \Rightarrow (iii)$ We may assume that
$\cl K =\cl{K(A)}= \sum_{j\in \cl J}\oplus_{c_o} \cl K(H_j)$
for some Hilbert spaces $H_j$. Let
$H_c = \vee_{j\in\cl J} H_j$.

Define $\cl C = \cl A/\cl K$ 
and note that $\cl C$ is generated by its projections. The
morphism $\tau : \cl C\rightarrow M(\cl K)/\cl K$ can be obtained
in the following way: Let $c\in\cl C$ and $a\in\cl A$ be any lift of $c$.
Then $\tau(c)=\pi(a|_{H_c})$. Indeed, since $\cl K$ acts non-degenerately on
$H_c$, we have that $a|_{H_c}$ is an element of the multiplier
algebra $M(\cl K)$  and the map $a\to a|_{H_c}$ is a morphism
extending the inclusion $\cl K\to M(\cl K)$, hence coincides with $\gs$. 
Writing $a=b+z$ with $b\in \cl K$ and $z\in\cl{Z(A)}$,
we see that $\tau(c) = \pi((b+z)|_{H_c})= \pi(z|_{H_c})$.
Since $\cl{Z(A)}$ acts `componentwise' on $\cl{K(A)}$, it follows
that $z|_{H_c}$ belongs to the centre of $M(\cl K)$ and so
$\tau(c)$ belongs to the centre of $M(\cl K)/\cl K$. 

\medskip

$(iii) \Rightarrow (i)$ Suppose that
$$0\rightarrow\cl K\rightarrow\cl A\rightarrow\cl C\rightarrow 0$$
is an exact sequence, $\cl C$ is abelian and generated by its
projections, $\cl K =\sum_{j\in \cl J}\oplus_{c_o} \cl K(H_j)$
and the corresponding morphism
$\tau : \cl C\rightarrow M(\cl K)/\cl K$ takes values in
the centre of $M(\cl K)/\cl K$.  The
projection onto $H_j$ is denoted $z_j$.

Recall that up to strong isomorphism of extensions \cite[3.2.11 and 3.2.13]{wo}  
we may write
$$\cl A = \{(a,c)\in M(\cl K)\oplus\cl C : \pi(a) =
\tau(c)\}.$$

\emph{We first show that $\cl A$ has the positive angle property.}

Any projection $P\in \cl A$ has the form $P = (p,q)$ where $p\in
M(\cl K)$ and $q\in\cl C$ are projections. Since
$\pi(p)=\tau(q)\in\cl Z((M(\cl K)) / \cl K)$, 
writing $p=\oplus_{j\in\cl J} p_j$ with each $p_j\in M(\cl K_j)=\cl B(H_j)$,
we have that
$\oplus_{j\in\cl J} (p_jx_j - x_jp_j) \in \cl K$, for all
$\oplus_{j\in\cl J} x_j\in M(\cl K)$. Now if $0 < p_j < z_j$, 
there exists $x_j\in\cl B(H_j)$ of norm
one such that $x_j=(z_j-p_j)x_jp_j$ and then $\nor{p_jx_j - x_jp_j}=\nor{x_j}=1$. 
It follows that the set $\cl J_0 = \{j\in\cl J : 0 < p_j < z_j\}$ is
finite. For each $j\in\cl J_0$, we easily verify that $p_j$ commutes with $\cl B(H_j)$
modulo $\cl K(H_j)$; since the Calkin algebra has trivial centre
\cite[Theorem 2.9]{cal}, this shows that $p_j=\gl_jz_j+a_j$ where
$\lambda_j\in\bb{C}$ and $a_j\in\cl K_j$. Hence $p_j$ must have
either finite rank or finite co-rank.

Take two projections $P_1 = (p_1,q_1)$ and $P_2 = (p_2,q_2)$ in
$\cl A$. By the previous paragraph, $c(p_1,p_2) < 1$; since $\cl C$
is abelian, $c(q_1,q_2) = 0$. Thus $c(P_1,P_2) = \max\{c(p_1,p_2),
c(q_1,q_2)\}< 1$ as required.

\smallskip
\emph{We now show that $\cl A$ is generated by its projections.}

Let $\cl A'$ be the C*-subalgebra of $\cl A$ generated by the projections of $\cl A$. 
Clearly $\cl A'$ contains $\cl K$. Since $\cl C$  is generated by its projections,
it suffices to show that every projection in $\cl C$ 
lifts to a projection in $\cl A$, hence in $\cl A'$.  
Indeed, it 
will then follow that $\cl A'/\cl K=\cl C=\cl A/\cl K$ and therefore $\cl A=\cl A'$. 

Thus let $q\in\cl C$ be a projection. 
Noting that $\tau(q)\in M(\cl K)/\cl K$, let $a\in M(\cl K)$ be any selfadjoint lift
of $\tau(q)$. Recall that 
$M(\cl K) = \sum_{j\in \cl J}\oplus_{\ell^{\infty}} \cl B(H_j)\subseteq \cl B(H)$. 
If $p\in \cl B(H)$ is the spectral projection of $a$ corresponding 
to the complement of the interval $[-1/2, 1/2]$, Calkin \cite[Theorem 2.4]{cal}
proves that $a-p$ is a compact operator; 
since $M(\cl K)$ is actually
a von Neumann algebra containing $a$, it follows that $p\in M(\cl K)$. Thus  
$a-p$ is a compact operator and is in $M(\cl K)$, hence it must belong to $\cl K$. 
Now $\pi(p)=\pi(a)=\tau(q)$, and hence $(p,q)\in \cl A$; this is the required lift of $q$
to $\cl A$.     $\quad\Box$

\begin{remark}
The sum $\cl{K(A)+Z(A)}$ in the  Theorem need not be direct. In fact, $\cl{K(A)}$
need not be complemented in $\cl A$ (example: $\cl
A=\ell^\infty)$.\\
Also, $\cl{Z(A)}$ need not have a complement which is an ideal of
$\cl A$ (example: $M_2)$.
\end{remark}

D. Topping has shown in \cite{to} that if $\cl A$ is a von Neumann
algebra such that  $c(p,q) <1$ for each pair of projections
$p,q\in \cl A$ then $\cl A$ is the direct sum of a finite
dimensional von Neumann algebra and an abelian one.
As a corollary to Theorem \ref{thm} we obtain the following
generalisation of  Topping's result:

\begin{proposition}\label{aw}
If $\cl A$ is an AW*-algebra with the positive angle property 
 then $\cl A$ is the direct sum of a
finite dimensional AW*-algebra and an abelian one.
\end{proposition}

\proof We use the notation introduced in the beginning of this section.
Let $\cl J_1=\{j\in \cl J: \dim H_j>1\}$. For each $j\in
\cl J_1$, the projection $z_j$ majorises  at least two minimal
orthogonal projections $e_j,f_j$ in $\cl A$. The completeness of
the lattice $\cl P(\cl A)$ implies that $p=\vee e_j$ is in
$\cl{P(A)}$. If $\cl J_1$ is infinite, this contradicts
Proposition \ref{l_zeroa}. Hence, $\cl J_1$ is finite.

If some
$z_j$ majorises an infinite sequence $\{e_n\}$ of pairwise orthogonal elements,
the projection $p=\vee_ne_{2n}$ is in $\cl A$, which again contradicts
Proposition \ref{l_zeroa}.
Thus $H_j$ must be finite dimensional for each $j\in\cl J_1$.

Let $z=\sum_{j\in\cl J_1}z_j$.
Then $z\cl A$ is finite-dimensional, $z^\bot \cl A\subseteq \cl{Z(A)}$ and
$\cl A$ is the direct sum of these two subalgebras.
\prend

\section{Equivalence of the Directed Set and the Positive Angle properties}

\begin{theorem}\label{dsp}
A C*-algebra has the directed set property if and only 
if it has the positive angle property.

\end{theorem}

\proof If a C*-algebra has the directed set property, then, for each pair
$p,q$ of projections, the algebra  $C^*(p,q)$ is clearly
finite-dimensional, and hence $c(p,q)<1$.

For the converse, since every finite-dimensional C*-algebra 
is generated by finitely many projections, it suffices to prove

\begin{theorem}\label{52}
Let $\cl A$ be a C*-algebra with the positive angle property which is generated by finitely many projections. Then $\cl A$ is finite-dimensional.
\end{theorem}
\proof 
 Note that $\cl A$ must be unital.
We use the notation of Section 4; in particular, recall that the projections $z_j$
commute with $\cl A\subseteq \bh$ and their sum is the identity operator
on $H$.

Fix a finite generating set $\cl F\subseteq \cl{P(A)}$.
Let $\cl J_0\subseteq \cl J$ be the set of
$j\in\cl J$ for which there exists
$p\in\cl F$ such that $0\ne pz_j\ne z_j$. It follows from the proof
of Proposition \ref{gen} that $\cl J_0$ is finite. 

Let $z=\sum_{j \in \cl J_0}z_j$. For each $p\in\cl F$ and each 
$j\in\cl J_1=\cl J\setminus\cl J_0$, the projection $pz_j$ is either $0$ or $z_j$. 
Since   
$$ z^{\perp}p = \sum_{j\in\cl J} z^{\perp}pz_j=\sum_{j\in\cl J_1} z^{\perp}pz_j$$
we see that $z^\perp p$ is a (finite) sum of $z_j$'s, hence commutes with $\cl A$.
In particular, the set $\{z^{\perp}p: p\in\cl F\}$
is commutative, hence so is the algebra $z^{\perp}\cl A$. Being generated by a finite
number of projections, it is finite-dimensional.

Therefore it remains to show that the algebra $z\cl A$ is also finite dimensional. 
This follows from the  

\smallskip

\noindent\textit{Claim }  
For each $j\in \cl J$, the algebra $z_j\cl A $ is
finite-dimensional.

\smallskip

\proofcl  
Since $\cl Az_j$ is equal to either $\cl K(H_j)$ or $\cl K(H_j)+\bb Cz_j$
(by Theorem \ref{thm}), each projection
$pz_j\in\cl Az_j$ is either of finite rank or of finite co-rank.

 Thus, $z_j\cl A$ is contained in a C*-algebra generated by the unit 
$z_j$ and finitely many projections
$p_1,\dots,p_k$ of finite rank. If $H_j'$ is the Hilbert space generated
by $p_iH_j, i=1,\dots, k$,  then $H_j'$ is finite dimensional and $z_j\cl A\subseteq
\cl B(H_j') + \bb{C}z_j$, which is finite dimensional.
 \prend

\end{document}